\documentclass[12pt]{article}
\usepackage{amsmath}
\numberwithin{equation}{section}
\usepackage{amsfonts,amssymb,latexsym,amsbsy, bbm, theorem,enumerate,color}
\usepackage{amssymb}
\usepackage{graphicx, graphics}
\usepackage{float,color,fancybox,shapepar,setspace,hyperref}
\usepackage{subfigure}
\usepackage{pgf,tikz}
\usetikzlibrary{arrows}
\usepackage{authblk}

\makeatletter

\newcommand{\Rmnum}[1]{\expandafter\@slowromancap\romannumeral #1@}
\makeatother
\newtheorem{Main Theorem}{Main Theorem}

\newtheorem{Conjecture}{Conjecture}
\newtheorem{Theorem}{Theorem}

\newtheorem{Lemma}{Lemma}
\newtheorem{Claim}{Claim}

\newtheorem{Proposition}{Proposition}
\newtheorem{definition}{Definition}[section]
\def\square{\hbox{\vrule height8pt depth0pt
\vbox{\hrule width7.2pt\vskip7.2pt\hrule width7.2pt}\vrule
height8pt depth0pt}\smallskip}

\def\pf{\medskip\noindent {\emph{\bf Proof}.}~~}

\newcommand{\ex}{{\rm  ex}}
\title{ Tur\'an problems  for suspension of a balanced tree}

\author[1]{Xiutao Zhu}

\author[2]{Xiaolin Wang}
\author[3]{Yanbo Zhang}
\author[4]{Fangfang Zhang\footnote{Corresponding author:fangfangzh@nufe.edu.cn}}

\date{}

\affil[1]{School of mathematics, Nanjing University of Aeronautics and Astronautics, Nanjing, 210016, P.R. China.}
\affil[2]{School of Mathematics and Statistics, Fuzhou University, Fuzhou 350108, P.R. China }
\affil[3]{School of Mathematical Sciences, Hebei Normal University, Shijiazhuang 050024, P.R. China }
\affil[4]{School of Applied Mathematics, Nanjing University of Finance and Economics, Nanjing 210023, China}

\begin{document}

\maketitle
\begin{abstract}
The Tur\'an number $\ex(n,H)$ is the maximum number of edges that an $n$-vertex $H$-free graph can have. The suspension $\widehat{H}$ is obtained from $H$ by adding a new vertex which is adjacent to all vertices of $H$ and a tree is balanced if the sizes of its two color classes differ at most $1$. In this paper, we obtain a sharp bound of $\ex(n,\widehat{T})$  when $n\ge 4(4k)^6$ based on the Erd\H{o}s-S\'os conjecture. We  also show the bound is sharp for infinitely many $n$ and characterize all extremal graphs. In particular, if $T$ satisfies some conditions such as $T$ contains a matching covering all vertices in one color class, then the bound is sharp for all $n$. This is a new class of graphs whose decomposition family does not contain a linear forest but we still can determine its Tur\'an number.
\vskip 2mm

\end{abstract}

\section{Introduction}
For a given family $\mathcal{H}$ of graphs, we say a graph $G$ is $\mathcal{H}$-free if it does not contain a subgraph isomorphic to any member in $\mathcal{H}$.  The Tur\'an number $\ex(n,\mathcal{H})$ is the largest number of edges in an $n$-vertex  $\mathcal{H}$-free graph. Usually, an $n$-vertex $\mathcal{H}$-free graph with $\ex(n,\mathcal{H})$ edges is called an extremal graph for $\ex(n,\mathcal{H})$. If $\mathcal{H}=\{H\}$, we denote it by $\ex(n,H)$. The studying of the Tur\'an number and the extremal graphs was initiated by Tur\'an \cite{turan} who determined the Tur\'an number $\ex(n,K_r)$ for any clique $K_r$. Let $T_r(n)$ denote the balanced complete $r$-partite graph on $n$ vertices where the difference in size between any two partite is at most $1$. Tur\'an proved the following result.

\begin{Theorem}(Tur\'an \cite{turan})\label{Turan}
For all $n$, $\ex(n,K_{r+1})=e(T_r(n))$ and $T_r(n)$ is the unique extremal graph.
\end{Theorem}

Since then, it has been an active problem in extremal graph theory and a highly significant result is Erd\H{o}s-Stone-Simonovits Theorem\cite{Er,stone}. For a general graph $H$, let $\chi(H)$ denote the chromatic number  of $H$.  Erd\H{o}s, Stone and Simonovits obtained the following equation:
\[\ex(n,H)=\left(1-\frac{1}{\chi(H)-1}\right)\binom{n}{2}+o(n^2).\]
Although the above equation shows us the asymptotic values for all non-bipartite graphs, we possess extremely limited information when $H$ is bipartite. Even the order of the magnitude of $\ex(n,H)$ is difficult to determine. Some studies about the magnitude of $\ex(n,H)$ can refer to \cite{Bukh,Jiang,Qiu,kang} and some other partial results can refer to \cite{Alon,Bondy,Kettle,Gallai,Furedi,Gao,Janzer,Kovari} and a survey in \cite{Furedi2}.

For the non-bipartite graphs, although Erd\H{o}s-Stone-Simonovits Theorem already tells the  asymptotic value of $\ex(n,H)$, it remains a challenge to determine the precise value of $\ex(n,H)$. However, for non-bipartite graphs, we have a potent tool, called decomposition family defined by Simionovits \cite{Simonovits2}, to assist us in handling such kind of problems. For a graph $H$ with $\chi(H)=r+1$, the decomposition family of $H$ is defined as following.
\begin{definition}
The decomposition family $\mathcal{M}(H)$ is the set of minimal graphs $M$ such that if a copy of $M$ is embedded into one partite set of $T_r(n)$, then the resulting graph contains $H$ as a subgraph.
\end{definition}
Simonovits \cite{Simonovits2} proved that if $\mathcal{M}(H)$ contains a linear forest, then the extremal graphs of $\ex(n,H)$ have very symmetrical structure. Up to now, most of the  Tur\'an numbers which can be calculated exactly rely on the condition that $\mathcal{M}(H)$ contains a linear forest. For example, the graph containing a color-critical edge \cite{Simonovits}, the edge blow-up graphs \cite{H.Liu,Wang,Yuan}, the odd-ballooning of a graph \cite{Hou,Peng,Yuan3,Zhu2,Zhu}, the pth-power of paths \cite{Xiao,Yuan2}, the graphs from geometric shapes, such as the icosahedron $I^{12}$ \cite{Simonovits3}, the dodecahedron $D^{20}$ \cite{Simonovits2} and the odd prism $C_{2k+1}\Box K_2$ \cite{He}. But if $\mathcal{M}(H)$ does not contain a linear forest, only a few  special Tur\'an numbers are known and we will introduce them later.

In this paper, we will discuss the Tur\'an number of the suspension of a graph. For a fixed graph $H$, the suspension of $H$, denoted by $\widehat{H}$, is the graph obtained from $H$ by adding a new vertex adjacent to all vertices of $H$. Consider it in this manner, the clique $K_{r+1}$ can be regarded as $\widehat{K}_r$. Inspired by a problem in connection with a geometrical problem, Erd\H{o}s posed the problem of $\ex(n, \widehat{T_2(6)})$. Later, this problem was solved by Simonovits\cite{Simonovits} in the general case of $\ex(n, \widehat{T_r(rd)})$. To show Simonovits's result, we define two families of graphs firstly. The graph $U_n$ is obtained by embedding a maximal $(d-1)$-regular graph\footnote{some one vertex has degree $d-2$ if such $(d-1)$-regular graph does not exist} into each partite of $T_r(n)$.  Then $\mathcal{U}_n$ consists of all such kind of $U_n$. If each partite of $U_n$ contains no triangles after embedding $(d-1)$-regular graphs, then let $U_n\in \mathcal{U}_n^*$. Obviously, $\mathcal{U}_n$ and $\mathcal{U}_n^*$ are not empty.

\begin{Theorem}(Simonovits \cite{Simonovits})\label{Simonovits}
Let $n$ be large enough. Then all graphs in $\mathcal{U}_n^*$ are extremal graphs of $\ex(n,\widehat{T_r(rd)})$ and all extremal graphs of $\ex(n,\widehat{T_r(rd)})$ must be in $\mathcal{U}_n$.
\end{Theorem}

Actually, Tur\'an Theorem \ref{Turan} is a special case of Simonovits's Theorem \ref{Simonovits}, since $K_{r+1}=\widehat{T_r(r)}$. These two theorems are the cases that $H$ is dense. On the other hand, the results on $\ex(n,\widehat{H})$ for $H$ being sparse are also varied. The first result in this case was obtained by Erd\H{o}s, F\"{u}redi, Gould and Gunderson \cite{friendship}. They determined the Tur\'an number $\ex(n,\widehat{kK_2})$ and the extremal graphs, where $kG$ denotes the disjoint union of $k$ copies of $G$. Later, this result was extended to $\ex(n,\widehat{kK_r})$ by Chen, Gould, Pfender and Wei \cite{Chen}. Recently, Hou, Li and Zeng \cite{HouJ} obtain a stronger result in this problem.   They obtained the extremal graph of $\ex(n,\widehat{H})$ when $H=H_1\cup\cdots\cup H_t$, where each $H_i$ is k-chromatic and contains a color-critical edge.   When $H=C_{2k+1}$ is an odd cycle($\widehat{H}$ is a wheel), an early theorem from \cite{Simonovits} shows $T_3(n)$ is the unique extremal graph for $\ex(n,\widehat{C}_{2k+1})$ when $n$ is large enough. For the even case $H=C_{2k}$, Yuan \cite{Yuan3} determined the value of $\ex(n,\widehat{C}_{2k})$ and all  extremal graphs. Furthermore, Yuan suggested studying the Tur\'an number $\ex(n,\widehat{P_k})$ at the end of \cite{Yuan3}.

In this paper, we deal with this problem in a general form. More precisely, we study the Tur\'an problem  $\ex(n,\widehat{T})$ for any balanced tree $T$. Here a balanced tree means the sizes of its two color classes differ at most one. So the path $P_k$ is a special balanced tree. And as previously discussed, it is not hard to prove that $\mathcal{M}(H)$ does not contain a linear forest when $H=T_r(rd)$ or $H=C_{2k}$. The graph $\widehat{T}$ is a new class whose decomposition family does not contain a linear forest.

 When studying the extremal problems about trees, a conjecture that has to be mentioned is the notoriously difficult Erd\H{o}s-S\'os Conjecture.
\begin{Conjecture}(Erd\H{o}s-S\'os Conjecture)
For any tree $T$, $\ex(n,T)\le \frac{|T|-2}{2}n$.
\end{Conjecture}
Here, $|T|$ denotes the number of vertices in $T$.  Since the conjecture was first proposed in 1962, a lot of efforts have been made, but it has not been solved. At present, only a few special cases of Erd\H{o}s-S\'os Conjecture were proved, such as the trees of diameter at most four \cite{Mclennan}, the trees containing a vertex which is the parent of at least $\frac{|T|-1}{2}$ leaves \cite{Sidorenko}. And in 2008, a proof of the conjecture was announced for very large trees by Ajtai, Koml\'os, Simonovits and Szemer\'edi. Some other results can refer to \cite{Balasubramanian,Brandt,Sacle}.

Let $$f(n,k)=\max\left\{n_0n_1+\left\lfloor\frac{(k-1)n_0}{2}\right\rfloor:n_0+n_1=n\right\}.$$
Based on the Erd\H{o}s-S\'os Conjecture, we have
\begin{Theorem}\label{main}
Let $T$ be a balanced tree of size $2k$ or $2k+1$ and Erd\H{o}s-S\'os Conjecture holds for all of its subtrees. When $n\ge 4(4k)^6$, we have
\begin{align}\label{eq1}
\ex(n,\widehat{T})\le f(n,k).
\end{align}
Moreover, the equality holds  for infinitely many $n$.
\end{Theorem}

The other parts of this paper are organized as follows. In Section 2, we will construct some graphs which make the equality in above theorem holds for infinitely many $n$'s. Furthermore, if the tree $T$ satisfies some conditions, then we will show that the equality holds for every $n$. Section 3 contains some preliminaries. We prove the  sharp upper bound and characterize the extremal graphs in Section 4.

\section{Some notation and the extremal graphs}
In this section, we will construct some $\widehat{T}$-free graph which makes the equality in Theorem \ref{main} holds for infinitely many $n$'s. Moreover, if the balanced tree $T$ satisfies some conditions, we will show the equality in Theorem \ref{main} holds for all $n$.

The idea for our constructions is inspired by the decomposition family. By the definition of the decomposition family, we can get the following simple proposition.
\begin{Proposition}\label{pro1}
If we embed an  $\mathcal{M}(H)$-free graph into one partite set of a large complete $r$-partite graph, then the resulting graph is $H$-free.
\end{Proposition}

To characterize the decomposition family of $\mathcal{M}(\widehat{T})$,  we first give some notation which would also be used in some other places. For a subset $A$ of $V(G)$, $G[A]$ denotes the subgraph of $G$ induced by the set $A$ and $G-A$ is induced by $V(G)\setminus A$. Let $T$ be a balanced tree. Say the small color class  of $T$ has size $k$. A covering $S$ of $T$ is a subset of $V(T)$ such that $V(T-S)$ is independent.   Let $\mathcal{S}=\{T[S]:~|S|\le k\}$ and $\widehat{\mathcal{S}}=\{\widehat{T[S]}: |S|\le k\}$. Obviously, the star $K_{1,k}$ is in $\widehat{\mathcal{S}}$, because $T$ has an independent covering of size $k$. However, any other member in $\widehat{\mathcal{S}}$ contains triangles since $T$ has no independent covering of size less than $k$. It is not hard to get the following.

\begin{Proposition}\label{pro2}
For any balanced tree $T$, $\mathcal{M}(\widehat{T})=\widehat{\mathcal{S}}\cup \{T\}$.
\end{Proposition}

By a basic calculation, we have
\begin{Proposition}\label{calculation}
The function $f(n,k)$ attains the maximum at
\begin{align}\label{eq2}
n_0=\left\{
\begin{array}{ll}
\frac{2n+k-1}{4} &~2n+k\equiv 1~\text{mod}~4 \\
\frac{2n+k-2}{4} &~2n+k\equiv 2~\text{mod}~4\\
\frac{2n+k-3}{4},\frac{2n+k+1}{4} &~2n+k\equiv 3~\text{mod}~4\\
\frac{2n+k}{4} &~2n+k\equiv 0~\text{mod}~4
\end{array}
\right.
\end{align}
\end{Proposition}

With the assistance of above propositions, we have,
\begin{Proposition}\label{proposition4}
Let  $T$ be a balanced tree on $2k$ or $2k+1$ vertices. If $2(k-1)|\lfloor\frac{2n+k}{4}\rfloor$, then equality in (\ref{eq1}) holds.
\end{Proposition}
\pf First let $K=K_{n_0,n_1}$ be a complete bipartite with two partite $A$ and $B$, where $|A|=n_0=\lfloor\frac{2n+k}{4}\rfloor$. Since $2(k-1)|\lfloor\frac{2n+k}{4}\rfloor$, we can embed as many as possible disjoint copies of $K_{k-1,k-1}$ into the partite $A$ and let $K^+$ denote the resulting graph. Obviously, $K^+[A]$ is a $(k-1)$-regular subgraph and so $K^+$  has
$$n_0n_1+\left\lfloor\frac{(k-1)n_0}{2}\right\rfloor$$ edges. On the other hand, $f(n,k)$ attains maximum at $n_0=\lfloor\frac{2n+k}{4}\rfloor$  by Proposition \ref{calculation}, so $e(K^+)=f(n,k)$, which is just the right hand of the inequality in Theorem \ref{main}.

We still need to prove $K^+$ is $\widehat{T}$-free. By proposition \ref{pro1} and \ref{pro2}, we just need to prove $K^+[A]$ is $\widehat{\mathcal{S}}\cup \{T\}$-free. This is evident since every member in $\widehat{\mathcal{S}}$ except $K_{1,k}$ contains triangles and $K^+[A]$ is $(k-1)$-regular but $\{T,K_3\}$-free. We are done. $\hfill\blacksquare$

\begin{Proposition}\label{proposition5}
Let  $T$ be a balanced tree on $2k+1$ vertices. If $2|\lfloor\frac{2n+k}{4}\rfloor$, then equality in (\ref{eq1}) holds.
\end{Proposition}
\pf The construction is a little similar to the above. Let $K=K_{n_0,n_1}$ be a complete bipartite with two partite $A$ and $B$, where $|A|=n_0=\lfloor\frac{2n+k}{4}\rfloor$. Firstly, embed as many as possible disjoint copies of $K_{k-1,k-1}$ into the partite $A$. Since $n_0$ is even, it is possible that there are $2t\le 2(k-2)$ remaining vertices  in $A$. Then we replace $t$ copies of $K_{k-1,k-1}$ in $A$ by $t$ copies of $K_{k,k}-kK_2$. Let $K^+$ denote the resulting graph. We have $K^+[A]$ is  $(k-1)$-regular and hence $e(K^+)=f(n,k)$.

 The graph $K^+$ is also $\widehat{\mathcal{S}}\cup \{T\}$-free since $K^+[A]$ is $(k-1)$-regular and each component is bipartite with at most $2k$ vertices. We are done. $\hfill\blacksquare$

\begin{Proposition}\label{proposition6}
Let  $T$ be a balanced tree on $2k$ or $2k+1$ vertices. If $T$ contains a matching of size $k$, then equality (\ref{eq1}) holds for every $n$.
\end{Proposition}
To prove this proposition, let $\beta(G)$ denote the size of the minimum covering of $G$ and  we need the following lemma.
\begin{Lemma}(K\"{o}nig \cite{Konig})\label{konig}
Let $G$ be a bipartite graph, then $\beta(G)=\nu (G)$.
\end{Lemma}
\pf When $T$ contains a matching of size $k$, then by Lemma \ref{konig}, $T$ contains no covering of size less than $k$. That is to say $\widehat{S}=\{K_{1,k}\}$ and $\mathcal{M}(\widehat{T})=\{K_{1,k}, T\}$. Now let $K=K_{n_0,n_1}$ be a complete bipartite with two partite $A$ and $B$, where $|A|=n_0=\lfloor\frac{2n+k}{4}\rfloor$. Since there exists a $(k-1)$-regular graph on $b$ vertices for any integer $b(\ge k)$, we can embed a $(k-1)$-regular graph into the partite $A$ of $K$ such that each component of the $(k-1)$-regular contains at most $2k-1$ vertices. The resulting graph has $f(n,k)$ edges and is $\widehat{T}$-free since the partite $A$ is $\{K_{1,k}, T\}$-free. $\hfill\blacksquare$

\section{Preliminaries}
In this section, we  do some preparatory work. Although the Erd\H{o}s-S\'os Conjecture is not resolved, we still have a weaker result which can guarantee the existence of a tree.

\begin{Lemma}\label{tree}
For any tree $T$ and any graph $G$ with $\delta(G)\ge |T|-1$, we can find a copy of $T$ in $G$.
\end{Lemma}

\begin{Lemma}\label{forest}
Let $k\ge k'\ge 0$ and $G$ be a graph on $n$ vertices with $\Delta(G)\le k-1$ and $e(G)\ge \frac{k'n}{2}+3k^3$. When $n\ge 3k^3$ and Erd\H{o}s-S\'os Conjecture is true, then
$$T_1\cup \cdots \cup T_N\subseteq G,$$
where $T_i$ is any tree of size at most $k'+2$ and $N\le 2k$.
\end{Lemma}
\pf We prove this lemma by induction on $N$. If $N=1$, then we can find any tree $T_1$ in $G$ by  Erd\H{o}s-S\'os Conjecture. Suppose we have found $ T_1\cup \cdots \cup T_{N-1}$ in $G$, then by $\Delta(G)\le k-1$
\[e(G-\cup_{i=1}^{N-1}T_i)\ge \frac{k'n}{2}+3k^3-(N-1)(k'+2)(k-1)> \frac{k'n}{2}.\]
Again by Erd\H{o}s-S\'os Conjecture, we can find any tree $T_N$ of size at most $ k'+2$ in $G-\cup_{i=1}^{N-1}T_i$. Together with $T_1\cup \cdots \cup T_{N-1}$, we find a copy of $T_1\cup \cdots \cup T_{N}$ in $G$. $\hfill\blacksquare$

\begin{Lemma}\label{decomposition}
Let $T$ be a balanced tree on $2k$ or $2k+1$ vertices. For any $a\in [k]$, there exists an independent set $I_a$ of size at most $a$ in $T$ such that
$$T-I_a=T_1\cup \cdots \cup T_N$$
with $|T_i|\le k-a+1$.
\end{Lemma}
\pf We prove the lemma by induction on $(k,a)$. First we consider two base cases $(2,a)$ and $(k,1)$.

For the case $k=2$, the structure of $T$ is as following.

\begin{figure}[h]
		\centering
		\begin{tikzpicture}
			\filldraw
			(0,1.5) circle (0.07)
			(1,1.5) circle (0.07)
            (0,0.5) circle (0.07)
			(1,0.5) circle (0.07)
			
            (3,1.5) circle (0.07)
            (4,1.5) circle (0.07)
            (2.5,0.5) circle (0.07)
            (3.5,0.5) circle (0.07)
            (4.5,0.5) circle (0.07)

            (6.5,1.5) circle (0.07)
            (7.5,1.5) circle (0.07)
            (6,0.5) circle (0.07)
            (7,0.5) circle (0.07)
            (8,0.5) circle (0.07) ;

            \draw[black, thick] (0,0.5) --(0,1.5)--(1,0.5) --(1,1.5);

            \draw[black,thick] (2.5,0.5)--(3,1.5)--(3.5,0.5)--(4,1.5)--(4.5,0.5) ;	

             \draw[black,thick] (6,0.5)--(6.5,1.5)--(7,0.5)  (6.5,1.5)--(8,0.5)--(7.5,1.5) ;	
			
		\end{tikzpicture}\quad
		
		\label{fig3}
	\end{figure}

\noindent It is not hard to find an independent set $I_a$ that makes the lemma holds.

Next consider the base case $(k,1)$. For the  balanced tree $T$, we see it as a bipartite graph with two color classes $A,B$ and $|A|=k$ and $|B|\in \{k,k+1\}$. Either we can find a leaf $x$ in $A$ and a leaf $y$ in $B$, or a leaf $y$ in $B$ and a vertex $x$ of degree $2$ in $A$ where $x$ is adjacent to $y$. The second case happens when $|B|=k+1$ and there are no leaves in $A$.  For both cases, we delete $x,y$ and obtain a smaller balanced tree $T'=T-\{x,y\}$. By induction hypothesis, there exists a vertex $w$ in $T'$ such that $T'-w=T_1'\cup \cdots \cup T_N'$ with $|T_i'|\le (k-1)-1+1=k-1$ for all $i\le N$.
Let us  consider adding  the two vertices $\{x,y\}$ back. If $x$ and $y$ are added into different components $T_i',T_j'$, then
$$\max\{|T_i'\cup \{x\}|,|T_j'\cup \{y\}|\}\le k=k-1+1$$
and hence $T-w$ satisfies the conditions. If $x$ and $y$ are added into the same component, say $T_1'$, and $|T_1'|\le k-2$, then $|T_1'\cup \{x,y\}|\le k=k-1+1$ and $T-w$ satisfies the conditions, too. If $|T_1'|= k-1$, then
$$|w\cup T_2'\cup \cdots \cup T_N'|\le 2k+1-(k+1)=k.$$
Let $w'$ be the neighbor of $w$ in $T_1'$. Then $T-w'= (T_1'\cup\{x,y\}-w')\cup (w\cup T_2'\cup \cdots \cup T_N')$ satisfies the condition, since each component of $T_1'\cup\{x,y\}$ has at most $k$ vertices.

At last consider the general case $k\ge 3$ and $a\ge 2$. Since $T$ is a tree, we can find a vertex $x$ such that all neighbors of $x$ except one are leaves. Let $x_1,\ldots,x_t$ be the neighbors of $x$ which are leaves and $y$ be the neighbor which is not a leaf. Without loss of generality,  we may assume $x\in A$ and let $x'$ be a neighbor of $y$ in $A-\{x\}$. First we delete $x$ and $x_t$, the vertices $x_1,\ldots,x_{t-1}$ become isolated. After that we add all edges between $x'$ and $x_1,\ldots,x_{t-1}$ and get a smaller balanced tree $T'$ on $2k-2$ or $2k-1$ vertices.
By the hypothesis, we can find an independent set $I'_{a-1}$ in $T'$ such that
\[T'-I'_{a-1}=T_1'\cup \cdots \cup T_N'~~\text{with}~~|T'_i|\le k-a+1~~\text{for all}~i.\]
Now we remove all edges $x'x_1,\ldots,x'x_{t-1}$ in $T'$ and then $x_1,\ldots,x_{t-1}$ become isolated vertices again and we try to add the vertices $x$ and $x_t$ back. If $y\notin I'_{a-1}$, then we put $x$ into the independent set and let
$I_a=I'_{a-1}\cup \{x\}-\{x_1,\ldots,x_{t-1}\}$. At this time,
\[T-I_{a}=T_1\cup \cdots \cup T_N\cup x_1\cup\cdots \cup x_{t},\]
where each component $T_i$ is obtained from $T_i'$ by deleting all possible leaves $\{x_1,\ldots,x_{t-1}\}$ and each $T_i$ satisfies $|T_i|\le k-a+1$. That is $T-I_a$ satisfies the conditions. If $y\in I'_{a-1}$, then $x'\notin I'_{a-1}$ since $x'$ is adjacent to $y$ in both $T'$ and $T$. We may assume $x'$ is contained in the component $T_1'$. Then some vertices of $\{x_1,\cdots,x_{t-1}\}$ are contained in $T_1'$ and some others are in $I'_{a-1}$. But  at most $k-a$ of $\{x_1,\cdots,x_{t-1}\}$ are contained in $T_1'$ since $|T_1'|\le k-a+1$. Now we consider the independent set $I_a=I'_{a-1}\cup \{x_t\}$, each component of $T-I_a$ is either a star induced by $x$ and some of $\{x_1,\cdots,x_{t-1}\}$ contained in $T_1'$,   or the original component $T_i'$ for $i\ge2$,  or $T_1'$ deleting some leaves in $\{x_1,\cdots,x_{t-1}\}$. All these components are of size at most $k-a+1$. That is $T-I_a$ satisfies the condition. We are done. $\hfill\blacksquare$

\section{The proof of main result}
Now we begin to prove our main result. We will use a standard trick to reduce the problem to graphs with large minimum degree.
\begin{Lemma}\label{largedegree}
Let $G$ be an $n$-vertex $\widehat{T}$-free graph with $\delta(G)\ge f(n,k)-f(n-1,k)$ and $n\ge 4(4k)^3$. Then
\[e(G)\le \max\left\{n_0n_1+\left\lfloor\frac{(k-1)n_0}{2}\right\rfloor:n_0+n_1=n\right\}.\]
\end{Lemma}

We first show how Lemma \ref{largedegree} implies Theorem \ref{main} and we give the proof of the lemma afterwards.

\noindent\textbf{Proof of Theorem \ref{main} using Lemma \ref{largedegree}}. Let $G$ be a $\widehat{T}$-free graph with $e(G)\ge f(n,k)$. If $\delta(G)\ge f(n,k)-f(n-1,k)$, then there is nothing to prove by Lemma \ref{largedegree}. Hence we may assume there is a vertex $v_n$ in $G$ with degree less than $f(n,k)-f(n-1,k)$. We define a process as follows: Let $G_n=G$ and $G_{n-1}=G_n-v_n$. If there is a vertex $v_{n-1}$ in $G_{n-1}$ with degree less than $f(n-1,k)-f(n-2,k)$, then delete it and get $G_{n-2}=G_{n-1}-v_{n-1}$. Proceed with this process until we get $G_\ell$ with $\delta(G_{\ell})\ge f(\ell,k)-f(\ell-1,k)$.

By this process, we have
\[e(G_i)\ge e(G_{i+1})-f(i+1,k)+f(i,k)+1, ~~~i\ge \ell,\]
and hence
\begin{align*}
e(G_{\ell})\ge& e(G_n)-f(n,k)+f(\ell,k)+(n-\ell)\\
\ge& f(\ell,k) +(n-\ell) \ge \left\lfloor\frac{\ell^2}{4}\right\rfloor+(n-\ell).
\end{align*}
On the other hand,
$e(G_{\ell})\le \binom{\ell}{2}$ and we have $n>\ell\ge \sqrt{4n+1}-1\ge 256k^3$.

Now by Lemma \ref{largedegree}, we have $e(G_{\ell})\le f(\ell, k)$ and so
\[e(G_n)\le e(G_{\ell})+\sum_{i=\ell}^{n-1}(f(i+1,k)-f(i,k)-1)\le f(n,k)-(n-\ell),\]
 a contradiction. We are done. $\hfill\blacksquare$

\vskip 5mm
\noindent\textbf{Proof of Lemma \ref{largedegree}}. Let $G$ be a graph satisfying the conditions in Lemma \ref{largedegree}. Hence $\delta(G)\ge f(n,k)-f(n-1,k)$.  When $2n+k\equiv 1~\text{mod}~4$ and by Proposition \ref{calculation},  we have
\begin{align*}
f(n,k)-f(n-1,k)=&\frac{(2n+k-1)(2n-k+1)}{16}+\frac{(k-1)(2n+k-1)}{8}\\
&-\frac{(2n+k-1)(2n-k-3)}{16}+\frac{(k-1)(2n+k-1)}{8}\\
\ge& \frac{(2n+k-1)}{4}\ge \frac{n}{2}.
\end{align*}
For all other $n$, we can prove the above inequalities by similar calculation and so $\delta(G)\ge \frac{n}{2}$. Next we will give some claims to characterize the structure of $G$ first.

\begin{Claim}\label{claim1}
We have $\Delta(G)\le \frac{n}{2}+2k$.
\end{Claim}
\pf Let $v$ be a vertex with maximum degree in $G$ and $N(v)$ be its neighborhood. Then $G[N(v)]$ is $T$-free and we can deduce that there exists a vertex $y\in N(v)$ such that $|N(y)\cap N(v)|\le |T|-2$ by Lemma \ref{tree}. Hence
\[n-\Delta(G)=|V(G)\setminus N(v)|\ge |N(y)|-|N(y)\cap N(v)|\ge \frac{n}{2}-(|T|-2),\]
which shows $\Delta(G)\le \frac{n}{2}+2k$ since $|T|\in\{2k,2k+1\}$. $\hfill\square$

\begin{Claim}\label{claim2}
We can find a partition $V(G)=V_0\cup V_1$ such that $V_0, V_1$ are nonempty and
\[\Delta(G[V_i])\le k-1~~~ \text{for}~~ i=0,1.\]
\end{Claim}
\pf Let $\delta(G)=\delta$, $\Delta(G)=\Delta$ and let $v$ be a vertex of degree  $\Delta$ in $G$. Define the sets $V_0$ and $V_1$ such that
\[|N(x)\cap N(v)|\ge \frac{2k-1}{2k}\Delta ~~~ \text{for}~~x \in V_0  \]
and
\[|N(x)\setminus N(v)|\ge\frac{2k-1}{2k}(n-\Delta) ~~~ \text{for}~~x \in V_1.\]

Now we show $V_0\cup V_1$ is a partition of $V(G)$. Obviously, $V_0\cap V_1=\emptyset$. Otherwise there exists a vertex of degree at least $\frac{2k-1}{2k}n$, contradicting  Claim \ref{claim1}. Both $V_0$ and $V_1$ are nonempty. Indeed, $v\in V_0$ and the vertex $y\in N(v)$ with $|N(y)\cap N(v)|\le |T|-2$(defined in Claim \ref{claim1})belongs to $V_1$, because
\begin{align*}
|N(y)\setminus N(v)|&\ge |N(y)|-(|T|-2)\\
&\ge\delta-(2k-1)\\
&\ge \frac{n}{2}-(2k-1)\\
&\ge \frac{2k-1}{2k}(n-\Delta).
\end{align*}
The last inequality holds since $\Delta\ge \delta\ge \frac{n}{2}$ and $n\ge 256k^3$. Finally, we prove every vertex $z$ is in $V_0\cup V_1$. If $|N(z)\cap N(v)|\le \frac{\Delta}{4k}$, then $z$ belongs to $V_1$ since for $n\ge 256k^3$ we have
\begin{align*}
|N(z)\setminus N(v)|\ge\delta-\frac{\Delta}{4k}\ge \frac{n}{2}-\frac{1}{4k}(\frac{n}{2}+2k)\ge \frac{2k-1}{2k}(n-\Delta).
\end{align*}
Next we consider the case $|N(z)\cap N(v)|> \frac{\Delta}{4k}$. Since both $G[N(z)]$ and $G[N(v)]$ are $T$-free, then by Lemma \ref{tree} we have
\[e(G[N(z)]\le  (|T|-2)|N(z)|\le (|T|-2)\Delta \]
and
\[e(G[N(v)]\le  (|T|-2)|N(v)|\le (|T|-2)\Delta. \]
By this, there are at most $\frac{1}{8k}\Delta$ vertices in $N(z)\cap N(v)$ who have at least $32k^2$ neighbors in $N(z)$ and at most another $\frac{1}{8k}\Delta$ vertices in $N(z)\cap N(v)$ who have at least $32k^2$ neighbors in $N(v)$. Otherwise
$$\max\{e(G[N(z)],e(G[N(v)]\}\ge 2k\Delta\ge (|T|-2)\Delta,$$
a contradiction. This implies we can find a vertex $w \in N(z)\cap N(v)$ such that
\[|N(w)-N(v)\cup N(z)|\ge \delta-64k^2\]
and so $|N(v)\cup N(z)|\le n-\delta +64k^2$. Now we have
\begin{align*}
|N(v)\cap N(z)|&=|N(v)|+|N(z)|-|N(v)\cup N(z)| \\
&\ge \Delta+\delta- (n-\delta +64k^2) \\
&\ge \Delta-64k^2\ge \frac{2k-1}{2k}\Delta,~~~~(\text{since~}n\ge 256k^3+8k^2)
\end{align*}
and so $z\in V_0$.

Therefore, above proofs show $V_0\cup V_1$ is indeed a partition of $V(G)$. Next we prove $\Delta(G[V_i])\le k-1$ for $i=0,1$. Without loss of generality, suppose $x\in V_0$ has $k$ neighbors, saying $\{x_1,x_2,\ldots,x_k\}$, in $V_0$. Then, $x,x_1,\ldots, x_k$ have
\[\Delta-(k+1)\frac{\Delta}{2k}\ge k+1 ~~~~(\text{since}~n\ge 256k^3  )\]
common neighbors in $N(v)$ by the definition of $V_0$. Let $y_1,\ldots,y_{k+1}$ be $k+1$ common neighbors of  $x,x_1,\ldots, x_k$. Then we can find a copy of $\widehat{T}$ using the  vertices $\{x\} \cup \{x_1,\ldots, x_k\} \cup \{y_1,\ldots,y_{k+1}\}$, a contradiction. Hence we have $\Delta(G[V_i])\le k-1$ for $i=0,1$ and the proof is complete. $\hfill\square$

\begin{Claim}\label{claim3}
Let $V(G)=V_0\cup V_1$ be the partition in Claim \ref{claim2}. Then
$$\frac{n}{2}+(k-1)\ge |V_i|\ge \frac{n}{2}-(k-1) ~~~~\text{for}~~i=0,1$$
 and each vertex in $V_i$ has at least $\frac{n}{2}-(k-1)$ neighbors in $V_{1-i}$.
\end{Claim}
\pf This claim follows immediately  from the truth  $\delta(G)\ge \frac{n}{2}$ and $\Delta(G[V_i])\le k-1$. \\ $\hfill\square$

\vskip 5mm
Now, suppose $V_0\cup V_1$ is the partition in Claim \ref{claim2} with $|V_0|=n_0\ge |V_1|= n_1$. Let $G_i=G[V_i]$ for $i=0,1$ and $G[V_0,V_1]$ be the bipartite graph induced by $V_0$ and $V_1$. We define a family of graphs $\mathcal{U}(n,T)$: we embed a $(k-1)$-regular $\widehat{\mathcal{S}}\cup \{T\}$-free graph into the larger partite of $K_{n_0,n_1}$, where $n=n_0+n_1$. Note that by Proposition \ref{proposition4}, \ref{proposition5} and \ref{proposition6}, $\mathcal{U}(n,T)$ is not empty if $\lfloor\frac {2n+k}{4}\rfloor$ satisfies some conditions. Let $\Delta(G_1)=k_1\le k-1$. We will continue the proof by distinguishing three cases and show the extremal graphs must be in $\mathcal{U}(n,T)$ when the equality (\ref{eq1}) holds with two additional graphs.

\textbf{Case 1} $k_1=0$.

In this case, we know $G_1$ is an empty graph and so
\[e(G)\le n_0n_1+\left\lfloor\frac{(k-1)n_0}{2}\right\rfloor\le f(n,k).\]
The equalities hold if and only if $n_0$ is the value in equality (\ref{eq2}), $G[V_0,V_1]$ is a complete bipartite graph, $G_0$ is $(k-1)$-regular and $\widehat{\mathcal{S}}\cup \{T\}$-free by the definition of decomposition family and Proposition \ref{pro2}. That is to say, $G\in \mathcal{U}(n,T)$.

\textbf{Case 2} $k_1\ge 1$ and $\Delta(G_0)\le k-k_1-1$.

In this case, we have
\begin{align}\label{eq4.0}
e(G)\le \left\lfloor\frac{(k-k_1-1)n_0}{2}\right\rfloor+\left\lfloor\frac{k_1n_1}{2}\right\rfloor+n_0n_1\le n_0n_1+\left\lfloor\frac{(k-1)n_0}{2}\right\rfloor.
\end{align}
Eventually, we have $e(G)\le f(n,k)$. Equalities in (\ref{eq4.0}) hold if and only if $G[V_0,V_1]$ is a complete bipartite graph, $n_1=n_0$ and one of $\{n_0,k_1,k\}$ is even or $n_1=n_0-1, k_1=1$ and $n_1,k$ are both even.
If $e(G)=f(n,k)$ also holds, then  we have $n_0=\lfloor\frac{2n+k}{4}\rfloor$ by Proposition \ref{calculation}, which implies $k\in \{2,3\}$ when $n_0=n_1$ and $k\in \{2,4\}$ when $n_1=n_0-1$.

For the case $n_0=n_1$, one of $\{n_0,k_1,k\}$ is even and $k\in \{2,3\}$. If $k=2$ or $k=3$ and $k_1=2$, then $G_0$ is an empty graph and $G_1$ is $(k-1)$-regular and $\widehat{\mathcal{S}}\cup \{T\}$-free. Thus $G\in \mathcal{U}(n,T)$. When $k=3$ and $k_1=1$, $G_0$ and $G_1$ are both a matching on $n_0=n_1$ vertices. However, in such a graph $G$, we still can find a copy of $\widehat{T}$ if $T$ is the balanced tree showed in Fig 2. Indeed, there is a vertex $v$ from each tree in Fig 2 such that each component of $T-v$ is of order at most $2$ and so we can find a copy of $\widehat{T}$ using one edge in $G_1$ and many disjoint edges in $G_0$. That is to say, for all other balanced tree on $6$ or $7$ vertices, we have an additional extremal graph, which is obtained from  $K_{n_0,n_1}$ by embedding a maximal matching into each partite.

\begin{figure}[h]
		\centering
		\begin{tikzpicture}
			\filldraw
		    (1, 2) circle (0.07) 	

			(0,1.5) circle (0.07)	
			(1,1.5) circle (0.07)
            (2,1.5) circle (0.07)

           (0,1) circle (0.07)
           (1,1) circle (0.07)

            (4, 2) circle (0.07) 	

			(3,1.5) circle (0.07)	
			(4,1.5) circle (0.07)
            (5,1.5) circle (0.07)

           (3,1) circle (0.07)
           (4,1) circle (0.07)
           (5,1) circle (0.07)

            (7.5, 2) circle (0.07) 	

			(6,1.5) circle (0.07)	
			(7,1.5) circle (0.07)
            (8,1.5) circle (0.07)
            (9,1.5) circle (0.07)

           (8,1) circle (0.07)
           (7,1) circle (0.07)

            ;

            \draw[black, thick] (0,1)--(0,1.5) --(1,2) --(1,1.5)--(1,1);
            \draw[black, thick] (1,2) --(2,1.5);
			\draw[black, thick] (3,1)--(3,1.5) --(4,2) --(4,1.5)--(4,1);
            \draw[black, thick] (4,2) --(5,1.5)--(5,1);
	\draw[black, thick] (7,1)--(7,1.5) --(7.5,2) --(8,1.5)--(8,1);
    \draw[black, thick] (6,1.5) --(7.5,2) --(9,1.5);
      \node at (4,0) {Fig 2};
			
		\end{tikzpicture}
	\end{figure}

For the case $n_1=n_0-1, k_1=1$, $n_1$ is even and $k\in \{2,4\}$. If $k=2$, then $G_0$ is an empty graph and $G$ is obtained by embedding a maximal matching into the small part of $K_{n_0,n_1}$\footnote{This extremal graph is a little different from the graphs in $\mathcal{U}(n,T)$ because the matching is embedded into the smaller part}. When $k=4$, $G_0$ is a $2$-regular and $G_1$ is a matching on $n_1$ vertices. For all balanced trees $T$ on $8$ or $9$ vertices except the trees in Fig 3, we always can find a vertex $v$ such that  each component of $T-v$ is of order at most $3$, or we can find two nonadjacent vertices $\{u,v\}$ such that each component of $T-\{u,v\}$ is of order at most $2$. And so we can find a copy of $\widehat{T}$ in $G$ using an edge in $G_1$ and many disjoint $P_3$'s in $G_0$, or using a $P_3$ in $G_0$ and many disjoint edges in $G_1$. That is to say, for the balanced trees in Fig 3, we have an additional extremal graph, where $G_0$ is a 2-regular graph and $G_1$ is a matching. Nevertheless, one caveat exists. If $T$ is one of the first three trees in Fig 3, then $G_0$ should be triangle-free. Otherwise, there exists two adjacent vertices $\{u,v\}$ in the first three trees such that each component of $T-\{u,v\}$ is of order at most $2$ and we can find a copy of $\widehat{T}$ using a triangle in $G_0$ and many disjoint edges in $G_1$, a contradiction.

\vskip 5mm
\begin{figure}[h]
		\centering
		\begin{tikzpicture}
			\filldraw	
			
			(1,1.5) circle (0.07)
            (3,1.5) circle (0.07)
            (1,0.5) circle (0.07)
           (1.5,0.5) circle (0.07)
           (0.5,0.5) circle (0.07)
           (3,0.5) circle (0.07)
           (3.5,0.5) circle (0.07)
           (2.5,0.5) circle (0.07)

           (5,1.5) circle (0.07)
            (7,1.5) circle (0.07)
           (5,0.5) circle (0.07)
           (5.5,0.5) circle (0.07)
           (4.5,0.5) circle (0.07)
           (7,0.5) circle (0.07)
           (7.5,0.5) circle (0.07)
           (6.5,0.5) circle (0.07)
           (8,0.5) circle (0.07)

            (9.5,1.5) circle (0.07)
            (11.5,1.5) circle (0.07)
            (12.5,1.5) circle (0.07)
           (9.5,0.5) circle (0.07)
           (10,0.5) circle (0.07)
           (9,0.5) circle (0.07)
           (11.5,0.5) circle (0.07)
           (12,0.5) circle (0.07)
           (11,0.5) circle (0.07)

            (0.5,-1) circle (0.07)
            (1.5,-1) circle (0.07)
            (2.5,-1) circle (0.07)
            (3.5,-1) circle (0.07)
            (0,-2) circle (0.07)
            (0.5,-2) circle (0.07)
            (1.5,-2) circle (0.07)
            (2.5,-2) circle (0.07)
            (3.5,-2) circle (0.07)

            (5,-1) circle (0.07)
            (6,-1) circle (0.07)
            (7,-1) circle (0.07)
            (8,-1) circle (0.07)
            (4.5,-2) circle (0.07)
            (5,-2) circle (0.07)
            (6,-2) circle (0.07)
            (7,-2) circle (0.07)
            (8,-2) circle (0.07)

            (9.5,-1) circle (0.07)
            (10.5,-1) circle (0.07)
            (11.5,-1) circle (0.07)
            (12.5,-1) circle (0.07)
            (9.5,-2) circle (0.07)
            (10.5,-2) circle (0.07)
            (11.5,-2) circle (0.07)
            (12.5,-2) circle (0.07)
            (9,-2) circle (0.07)

            (0.5,-3.5) circle (0.07)
            (1.5,-3.5) circle (0.07)
            (2.5,-3.5) circle (0.07)
            (3.5,-3.5) circle (0.07)
            (0.5,-4.5) circle (0.07)
            (1.5,-4.5) circle (0.07)
            (2.5,-4.5) circle (0.07)
            (3.5,-4.5) circle (0.07)
            (4,-4.5) circle (0.07)

            (5,-3.5) circle (0.07)
            (6,-3.5) circle (0.07)
            (7,-3.5) circle (0.07)
            (8,-3.5) circle (0.07)
            (5,-4.5) circle (0.07)
            (6,-4.5) circle (0.07)
            (7,-4.5) circle (0.07)
            (8,-4.5) circle (0.07)
            (8.5,-4.5) circle (0.07)

             (9.5,-3.5) circle (0.07)
             (10,-3.5) circle (0.07)
            (11,-3.5) circle (0.07)
            (12,-3.5) circle (0.07)
             (10,-4.5) circle (0.07)
            (11,-4.5) circle (0.07)
            (11.5,-4.5) circle (0.07)
           (12,-4.5) circle (0.07)
            (12.5,-4.5) circle (0.07)
            ;

            \draw[black, thick] (1.0,0.5)--(1,1.5)--(3,1.5) --(3.0,0.5) (0.5,0.5)--(1,1.5)--(1.5,0.5) (2.5,0.5)--(3,1.5)--(3.5,0.5) ;

            \draw[black, thick] (5.0,0.5)--(5,1.5)--(7,1.5) --(7.0,0.5) (4.5,0.5)--(5,1.5)--(5.5,0.5) (6.5,0.5)--(7,1.5)--(7.5,0.5) (8,0.5)--(7,1.5);
            \draw[black, thick] (9.5,0.5)--(9.5,1.5)--(11.5,1.5) --(11.5,0.5) (9,0.5)--(9.5,1.5)--(10,0.5) (11,0.5)--(11.5,1.5)--(12,0.5)--(12.5,1.5);
            \draw[black, thick] (0.5,-1)--(0.,-2)   (0.5,-1)--(0.5,-2)  (0.5,-1)--(1.5,-2) --(1.5,-1)--(2.5,-2)--(2.5,-1)--(3.5,-2)--(3.5,-1);
      \draw[black, thick]  (6,-1)--(6,-2) (4.5,-2)--(5,-1)--(5,-2)  (5,-1)--(6,-2) --(7,-1)--(7,-2)--(8,-1)--(8,-2);
\draw[black, thick] (9,-2)--(9.5,-1)--(9.5,-2)--(10.5,-1)--(11.5,-2)--(11.5,-1)--(12.5,-2)--(12.5,-1)  (10.5,-1)--(10.5,-2);
\draw[black, thick] (0.5,-4.5)--(0.5,-3.5)--(1.5,-4.5)--(2.5,-3.5)--(3.5,-4.5)--(3.5,-3.5)--(4,-4.5) (1.5,-4.5)--(1.5,-3.5) (2.5,-3.5)--(2.5,-4.5);

\draw[black, thick] (5,-4.5)--(5,-3.5)--(6,-4.5)--(6,-3.5)--(7,-4.5)--(7,-3.5)--(8,-4.5)--(8,-3.5)--(8.5,-4.5);

\draw[black, thick] (9.5,-3.5)--(10,-4.5)--(11,-3.5)--(11.5,-4.5)--(12,-3.5)--(12.5,-4.5)  (10,-3.5)--(10,-4.5) (11,-3.5)--(11,-4.5) (12,-3.5)--(12,-4.5);

     \node at (6,-6) {Fig 3};
			
		\end{tikzpicture}\quad
		
		\label{fig3}
	\end{figure}

\textbf{Case 3} $k_1\ge 1$ and $\Delta(G_0)\ge k-k_1$.

We have $n_i\ge 128k^3$ for $i=0,1$  by Claim \ref{claim3}.
Let $y\in V_0$ be a vertex of degree at least $k-k_1$ in $G_0$ and $\{y_1,\ldots,y_{k-k_1}\}\subseteq N(y)\cap V_0$. We define
$$Y=N(y)\cap N(y_1)\cap \cdots\cap N(y_{k-k_1})\cap V_1$$
and by Claim \ref{claim3}, $y$ and each $y_i$ have at most $k-1$ non-neighbors in $V_1$, so
\begin{align}\label{eq4.1}
|Y|\ge n_1-(k-1)(k-k_1+1)\ge n_1-k^2\gg  3k^3.
\end{align}
By Lemma \ref{decomposition}, we can find an independent set $I_{k-k_1}$ in $T$ such that
$$T-I_{k-k_1}=T_1\cup \cdots\cup T_N~~\text{with}~~|T_i|\le k_1+1 ~~\text{and}~~N\le 2k.$$
If $e(G[Y])\ge \frac{k_1-1}{2}|Y|+3k^3$, then we can find a copy of $T_1\cup \cdots\cup T_N$ in $G[Y]$ by Lemma \ref{forest}. This copy of $T_1\cup \cdots\cup T_N$, together with the independent set $\{y_1,\ldots,y_{k-k_1}\}$, induce a copy of $T$ in the neighborhood of $y$, and so we find a copy of $\widehat{T}$, a contradiction. Hence, we have $e(G[Y])<\frac{k_1-1}{2}|Y|+3k^3$ and
\begin{align}\label{eq4.2}
e(G_1)< e(G[Y])+(n_1-|Y|)k_1\le \frac{k_1-1}{2}n_1+\frac{(k_1+1+6k)k^2}{2}.
\end{align}
The second inequality follows from the inequality (\ref{eq4.1}).

Now recall that $\Delta(G_1)=k_1$. Let $x\in V_1$ be the vertex of degree $k_1$ in $G_1$ and $N(x)\cap V_1=\{x_1,\ldots,x_{k_1}\}$. Analogously, we define
$$X=N(x)\cap N(x_1)\cap \cdots\cap N(x_{k_1})\cap V_0$$
and by Claim \ref{claim3}, $x$ and each $x_i$ have at most $2(k-1)$ non-neighbors in $V_0$, so
\begin{align}\label{eq4.3}
|X|\ge n_0-2(k-1)(k_1+1)\ge n_0-2k^2\gg  3k^3.
\end{align}
By Lemma \ref{decomposition}, we can find an independent set $I_{k_1}$ in $T$ such that
$$T-I_{k_1}=T_1\cup \cdots\cup T_{N'}~~\text{with}~~|T_i|\le k-k_1+1~~\text{and}~~N'\le 2k.$$
If $e(G[X])\ge \frac{k-k_1-1}{2}|X|+3k^3$, then we can find a copy of $T_1\cup \cdots\cup T_{N'}$ in $G[X]$ by Lemma \ref{forest}. This copy of $T_1\cup \cdots\cup T_{N'}$, together with the independent set $\{x_1,\ldots,x_{k_1}\}$, induce a copy of $T$ in the neighborhood of $x$, and so we find a copy of $\widehat{T}$, a contradiction.
Hence we have $e(G[X])< \frac{k-k_1-1}{2}|X|+3k^3$ and
\begin{align}\label{eq4.4}
e(G_0)< e(G[X])+(n_0-|X|)(k-1)\le \frac{k-k_1-1}{2}n_0+(4k+k_1+1)k^2.
\end{align}
The second inequality follows from the inequality (\ref{eq4.3}).

Now, combining the inequalities (\ref{eq4.2}) and (\ref{eq4.4}), we have
\begin{align*}
&e(G)\le e(G_1)+e(G_0)+n_0n_1\\
<&\frac{k_1-1}{2}n_1+\frac{(k_1+1+6k)k^2}{2}+\frac{k-k_1-1}{2}n_0+(4k+k_1+1)k^2+n_0n_1\\
\le& \frac{k-1}{2}n_0+n_0n_1-\frac{1}{2}n_1+9k^3\\
\le& n_0n_1+\left\lfloor\frac{(k-1)n_0}{2}\right\rfloor.
\end{align*}
The last inequality holds because $\frac{n_1}{2}\ge \frac{n}{2}-(k-1)\gg 9k^3$. We complete the proof.$\hfill\blacksquare$

\section{Concluding Remarks}
$\bullet$ In Section 2, we construct some $\widehat{T}$-free graphs to show that the equality (\ref{eq1}) holds for infinitely many values of $n$. In particular, if $T$ contains a matching that covers one color class, then the equality (\ref{eq1}) holds for every $n$. However, for the other trees, we do not know if the equality still holds for all $n$. If we want to prove the equality (\ref{eq1}) holds for every $n$, we need to find a $(k-1)$-regular and $\widehat{S}\cup \{T\}$-free graph on $n_0$ vertices. As pointed out, the other members of $\widehat{S}$ except $K_{1,k}$ contain triangles, so a method for constructing such a $(k-1)$-regular graph is to find a $\{K_3,T\}$-free $(k-1)$-regular graph. However, we can only find some specific $\{K_3,T\}$-free $(k-1)$-regular graph, such as the disjoint union of $K_{k-1,k-1}$ or $K_{k,k}-kK_2$ if $|T|=2k+1$. The challenge is to find such a $(k-1)$-regular graph of size $n_0$ for every $n_0$. We don't even know if such a graph exists.

By the way, if we revise a proof by Bushaw and Kettle (see Lemma 3.5 in  \cite{Kettle}), it is not hard to prove that for most of balanced trees $T$ without ``perfect matching" , the condition $(k-1)$-regular and $\widehat{S}\cup \{T\}$-free is equivalent to $(k-1)$-regular and $\{K_3,T\}$-free.

$\bullet$ In our result, we need the condition that $T$ is balanced. However, for any unbalanced tree whose smaller color class has size $k$, Theorem \ref{main} still holds if Lemma \ref{decomposition} holds for this tree. We don't even need to revise the proof. For example, Lemma \ref{decomposition} holds for any double star $S_{k,t}$ with $t\ge k$ and so $\ex(n,\widehat{S}_{k,t})\le f(n,k)$. Based on this, we conjecture that our result holds for any tree.

\begin{Conjecture}
Let $T$ be a tree such that the smaller color class has size $k$. Then for large $n$,
\[\ex(n,\widehat{T})\le f(n,k).\]
\end{Conjecture}

\section{Acknowledgment}

The research is supported by  National Natural Science Foundation of China under number 12101298, 12401447, 12401454.  The author Zhu is also supported by Basic Research Programma of Jiangsu Province(BK20241361), Jiangsu Funding Program for Excellent Postdoctoral Talent(2024ZB179), State-sponsored Postdoctoral Researcher program(GZB20240976).

\end{document}